\documentclass[12pt]{article}
\usepackage{amsfonts} 
\usepackage{theorem}  
\usepackage{amsmath}
\newtheorem{thm}{Theorem}[section]

\newtheorem{prop}[thm]{Proposition}

{\theorembodyfont{\normalfont\rmfamily}
\newtheorem{defn}[thm]{Definition}
\newtheorem{rem}[thm]{Remark}

}

\makeatletter
\newenvironment{introthm}[1]%
 {\renewcommand{\theththm}{\ref{#1}}\begin{ththm}}
 {\end{ththm}}
\newtheorem{ththm}{Theorem}
\makeatother
\makeatletter
\newenvironment{pf}[1][]%
 {\def\proof@temp{#1}\par\noindent
  \textsc{Proof}\ifx\proof@temp\@empty\else\ (#1)\fi\hspace{1em}}
 {\rule{12pt}{0pt}\hfill\rule{6pt}{6pt}\par\vspace{.4\baselineskip}}
\makeatother

\makeatletter
\def\operatorname#1{\mathop{\operator@font #1}\nolimits}%
\makeatother
\def\map#1 {\stackrel{#1}{\longrightarrow}}

\renewcommand{\H}{\mathbb{H}}

\newcommand{\R}{\mathbb{R}}

\newcommand{\gl}{\mathfrak{gl}}
\renewcommand{\sl}{\mathfrak{sl}}

\newcommand{\B}{\mathcal{B}}
\newcommand{\ca}{\mathcal{D}} 
\renewcommand{\H}{\mathcal{H}}
\newcommand{\J}{\mathcal{J}}

\newcommand{\T}{\mathcal{T}}
\newcommand{\V}{\mathcal{V}}
\newcommand{\W}{\mathcal{W}}

\newcommand{\End}{\operatorname{End}}

\newcommand{\Hom}{\operatorname{Hom}}
\newcommand{\Id}{\operatorname{Id}}

\newcommand{\Tr}{\operatorname{Tr}}

\newcommand{\rpn}{\R P^n}

\title{Affine connections with $W=0$}
\author{%
Francis Burstall\\[5pt]
\small\ttfamily F.E.Burstall@maths.bath.ac.uk\\[1pt]
\small Mathematical Sciences\\[-10pt]
\small University of Bath\\[-10pt]
\small Bath\ \ BA2\ 7AY\\[-10pt]
\small United Kingdom\\[20pt]
John Rawnsley\\[5pt]
\small\ttfamily J.Rawnsley@warwick.ac.uk\\[1pt]
\small Mathematics Institute\\[-10pt]
\small University of Warwick\\[-10pt]
\small Coventry\ \ CV4\ 7AL\\[-10pt]
\small United Kingdom\\[10pt]~
}

\date{26 February 2007}

\renewcommand{\baselinestretch}{1.3}

\begin{document}
\maketitle

\renewcommand{\baselinestretch}{1}
\vspace{.2in}
\begin{abstract}
If $\nabla$ is a torsionless connection on the tangent bundle of a
manifold $M$ the Weyl curvature $W^\nabla$ is the part of the curvature
in kernel of the Ricci contraction. We give a coordinate free proof of
Weyl's result that $W^\nabla$ vanishes if and only if $(M,\nabla)$ is
(locally) diffeomorphic to $\rpn$ with $\nabla$, when transported to
$\rpn$, in the projective class of $\nabla_{LC}$, the Levi-Civita
connection of the Fubini--Study metric on $\rpn$.

If $M$ is even dimensional and $J(M)$ denotes the bundle of all
endomorphisms $j$ of the tangent spaces of $M$, a connection $\nabla$
determines an almost complex structure $J^\nabla$ on $J(M)$ \cite{bib:OBR}.
We show that $J^\nabla$ is a projective invariant, that an integrable
$J^\nabla$ can be obtained from a torsionless connection and that we
must then have $W^\nabla=0$. We also show for torsionless connections
$\nabla$, $\nabla'$ that $J^\nabla = J^{\nabla'}$ if and only if 
$\nabla$ and $\nabla'$ are projectively equivalent.
\end{abstract}
\thispagestyle{empty}
\newpage
\renewcommand{\baselinestretch}{1.3}

\setcounter{page}{1}

\section{Introduction}

In Riemannian geometry the Ricci tensor splits into two irreducible
pieces under the orthogonal group, the scalar curvature which is its
trace and the traceless Ricci tensor. The part of the curvature tensor
lying in the kernel of the Ricci contraction is also irreducible and
known as the Weyl tensor. If the manifold is oriented, this
decomposition is still irreducible under the special orthogonal group
except in dimension $4$ where the Weyl tensor decomposes into its self-
and anti-self-dual parts. This decomposition was obtained by Singer and
Thorpe \cite{bib:SingThor}. It is shown in Eisenhart
\cite{bib:Eisenhart} that the Weyl tensor is the obstruction to the
Riemannian manifold being conformally flat.

The bundle of Hermitean structures on the tangent spaces of an even
dimensional Riemannian manifold carries a natural almost complex
structure whose integrability condition is the vanishing of the Weyl
tensor \cite{bib:BBO,bib:Dubois-Violette,bib:OBR}. In dimension 4, if
only Hermitean structures compatible with an orientation are used then
the integrability condition is the vanishing of the self-dual part of
the Weyl tensor which leads to the celebrated Riemannian analogue of
Penrose's twistor theory developed by Atiyah, Hitchin and Singer
\cite{bib:AHS}.

Other authors have considered decompositions of the curvature tensor
into irreducible components and corresponding integrability conditions
in a number of contexts, in particular for unitary
\cite{bib:OBR,bib:TricVanh}, quaternionic \cite{bib:Salamon} and
symplectic structures \cite{bib:Vaisman1,bib:Vaisman2}.

In the present paper we look at the case with the least restriction
on the structure group, the case of a linear connection on a manifold
and show that this leads naturally to projective geometry.

A torsionless connection $\nabla$ on the tangent bundle of an
$n$-dimensional manifold $M$ defines a \emph{projective structure}:
the family of connections sharing the same geodesics as $\nabla$.
The Weyl curvature tensor $W^\nabla$ of $\nabla$ is the part of the
curvature in kernel of the Ricci contraction.  When $n=2$,
$W^\nabla\equiv0$ but there is an analogue $C^\nabla$ of the
Cotton--York tensor of $3$-dimensional conformal geometry.  We give a
low-technology, coordinate-free proof of Weyl's theorem
\cite{bib:Weyl} that $W^\nabla$ and, for $n=2$, $C^\nabla$, are
projectively invariant and vanish if and only if $(M,\nabla)$ is
(locally) isomorphic to $(\rpn,\nabla_{LC})$ where $\nabla_{LC}$ is
the Levi-Civita connection.

If $J(M)$ denotes the bundle of all endomorphisms $j$ of the tangent
spaces of $M$, a connection $\nabla$ determines an almost complex
structure $J^\nabla$ on $M$. We show that $J^\nabla$ is a projective
invariant, that an integrable $J^\nabla$ can always be obtained from
a torsionless connection and, in that torsionless case, $J^\nabla$ is
integrable if and only if $W^\nabla=0$. In particular, when $n=2$,
$J^\nabla$ is always integrable even when $C^\nabla$ is non-zero.  We
also show that for torsionless connections $\nabla$ and $\nabla'$ we
have $J^\nabla = J^{\nabla'}$ if and only if $\nabla$ and $\nabla'$
are projectively equivalent.

The paper is structured as follows:

In section 2 we set the scene and establish notation.

In section 3 we describe the decompositions of torsion and curvature
tensors and show that the Weyl component of the curvature is
projectively invariant.

In section 4 we look at the twistor theory and show that the Weyl
component of the curvature is the obstruction to integrability of the
twistor almost complex structure.

\begin{introthm}{thm-proj-inv}
The almost complex structure $J^\nabla$ on $J(M)$ defined by a
connection $\nabla$ on $M$ only depends on the projective class of
$\nabla$. If $\nabla$ and $\nabla'$ both have zero torsion then
$J^\nabla = J^{\nabla'}$ if and only if $\nabla$ and $\nabla'$ are
projectively equivalent.
\end{introthm}

\begin{introthm}{thm-torzero}
If $J^\nabla$ is integrable then there is another connection
$\nabla'$ defining the same almost complex structure on $J(M)$ and with
zero torsion.
\end{introthm}

\begin{introthm}{thm-intgrbl}
Let $\nabla$ be a torsion-free connection then $J^\nabla$ is
integrable if and only if $W^\nabla = 0$.
\end{introthm}

In section 5 we prove
\begin{introthm}{thm-projflat}[\cite{bib:Weyl}]
Let $M$ be an $n$-dimensional manifold with tor\-sion-free connection
$\nabla$.  Suppose that either $n\geq3$ and $W^\nabla=0$ or $n=2$ and
$C^\nabla=0$.  Then there are local affine diffeomorphisms
between $M$ and $\rpn$ equipped with a connection in the projective
class of the symmetric connection.
\end{introthm}

\section{Preliminaries}

\subsection{Connections, curvature and torsion}
\label{sec:conn-curv-tors}

Let $M$ be a manifold and $\nabla$ a connection in $TM$. Its torsion
$T^\nabla$ and curvature $R^\nabla$ are given by
\[
\begin{array}{lcl}
T^\nabla(X,Y) &= &\nabla_XY-\nabla_YX-[X,Y],\\  R^\nabla(X,Y)Z &=
&\nabla_X\nabla_YZ - \nabla_Y\nabla_XZ -\nabla_{[X,Y]}Z.
\end{array}
\]
$T^\nabla$ is a $TM$-valued 2-form and $R^\nabla$ an $\End TM$-valued
$2$-form. As forms they have covariant exterior derivatives computed
using $\nabla$.
It is easy to check that we have the First Bianchi Identity
\begin{equation}\label{Bianchi1}
(d^\nabla T^\nabla)(X,Y,Z)= R^\nabla(X,Y)Z + R^\nabla(Y,Z)X +
R^\nabla(Z,X)Y
\end{equation}
and the Second Bianchi Identity
\begin{eqnarray}\label{Bianchi2}
0 &=& (d^\nabla R^\nabla)(X,Y,Z)\nonumber\\
&=& (\nabla_XR^\nabla)(Y,Z) + (\nabla_YR^\nabla)(Z,X) +
(\nabla_ZR^\nabla)(X,Y)\nonumber\\
&&\qquad\mbox{} + R(T(X,Y),Z) + R(T(Y,Z),X) + R(T(Z,X),Y).
\end{eqnarray}

The \textit{Ricci curvature}, $r^\nabla$, of a linear connection is
given by
\[
r^\nabla(X,Y) = \Tr\left( Z \mapsto R^\nabla(X,Z)Y\right).
\]
There is a second trace we could take
\[
s^\nabla(X,Y) = \Tr(R^\nabla(X,Y))
\]
which gives a 2-form. If $X_i$ is a local frame field and $\alpha^i$ the
dual frame field so that $\alpha^i(X_j) = \delta^i_j$ then, using the
First Bianchi Identity,
\begin{align*}
s^\nabla(X,Y) &= \sum_i \alpha^i(R(X,Y)X_i)\\
&= \sum_i \alpha^i\biggl((d^\nabla T^\nabla)(X,Y,X_i) - R(Y,X_i)X -
R(X_i,X)Y\biggr)\\
&= r^\nabla(X,Y) - r^\nabla(Y,X) + \sum_i \alpha^i\biggl((d^\nabla
T^\nabla)(X,Y,X_i)\biggr).
\end{align*}
In particular, the Ricci tensor is symmetric when $s^\nabla=0$ and
$\nabla$ is torsion-free, but not in general.

When $\nabla$ is torsion free, the second trace is determined by the 
antisymmetric part of the Ricci tensor.

\subsection{A bundle of Lie algebras}
\label{sec:bundle-lie-algebras}

The bundle $TM\oplus\End TM\oplus T^*M$ carries the
structure of a bundle of Lie algebras that will be useful to us.  For
this, declare $TM$ and $T^*M$ to be abelian subalgebras, give $\End
TM$ the usual commutator bracket and then, for $(X,A,\alpha)\in
TM\oplus\End TM\oplus T^*M$, set
\begin{align*}
[A,X] &= AX\\
[ \alpha ,A] &= \alpha\circ A
\end{align*}
and define $[X,\alpha]\in\End TM$ by
\[
[X,\alpha]Y=\alpha(X)Y+\alpha(Y)X.
\]
It is straightforward to check that this bracket does indeed satisfy
the Bianchi identity but a more conceptual explanation is available:
fix a line bundle $\Lambda$ and set $V=\Lambda\oplus TM\Lambda$ where
here and below we use juxtaposition to denote tensor product with a
line bundle.  Contemplate the bundle $\sl(V)$ of trace-free
endomorphisms of $V$.  This bundle of Lie algebras decomposes:
\[
\sl(V)=\Hom(\Lambda,TM\Lambda)\oplus
\sl(\End(\Lambda)\oplus\End(TM\Lambda))\oplus
\Hom(TM\Lambda,\Lambda).
\]
The first and last summands are canonically isomorphic to $TM$ and
$T^*M$ respectively, while the adjoint action provides an isomorphism
\[
\sl(\End(\Lambda)\oplus\End(TM\Lambda))\cong\End(\Hom(\Lambda,TM\Lambda))
\cong\End TM.
\]
Putting all this together, we arrive at a bundle isomorphism
\[
\sl(V)\cong TM\oplus\End TM\oplus T^*M
\]
which is readily shown to be an isomorphism of Lie algebras.

\subsection{$G$-structures and Representations}

Let $G$ be a Lie group and $M$ a manifold of dimension $n$ then a
$G$-structure on $M$ is a principle $G$-bundle $\pi\colon P \to M$
together with an $n$-dimensional representation $V$ of $G$ such that
$TM$ is isomorphic to the associated bundle $P\times_G V$. More
precisely,  we consider the $V$-frame bundle $Fr(M)$ consisting all the
linear isomorphisms $b \colon V \to T_xM$ of $V$ with the tangent spaces
of $M$ and the obvious projection map $\pi_M$, then a $G$-structure is a
morphism $P \to Fr(M)$ covering the homomorphism $G \to GL(V)$ given by
the representation. The isomorphism of $P\times_G V$ with $TM$ is then
induced by the identification $TM = Fr(M)\times_{GL(V)} V$.

This leads to further identifications of space of tensors on $M$ with
associated bundles of $P$. For instance, torsion tensors are sections of
the bundle associated to the representation $\Lambda^2V^*\otimes V$ and
curvature tensors to $\Lambda^2V^*\otimes \gl(V) = \Lambda^2V^*\otimes
V^*\otimes V$. The $G$-structure in this case is the standard one with
$G=GL(V)$ and $P=Fr(M)$. We adopt the more general language simply
because the questions we examine here make sense in the more general
context.

Even when the initial representation of $G$ on $V$ is irreducible, the
representations on tensor spaces may not be. For example $\Lambda^k V$
and $\Lambda^k V^*$ are irreducible under $GL(V)$ but $\Lambda^2 V^*
\otimes V$ is not. If $W$ is some representation of $G$ and $W = \bigoplus_k
W_k$ a decomposition into irreducible subrepresentations (we assume here
that we are dealing with groups for which such a decomposition is always
possible and unique up to multiplicities) then this splitting induces a
splitting of the corresponding associated bundles and we can project
$P\times_G W$ into the subbundles $P\times_G W_k$. We refer to the
projections of a section of $P\times_G W$ into the $P\times_G W_k$ as
its \textit{irreducible components}.

In the next section we determine the irreducible components of tensors
of torsion and curvature type for $G=GL(V)$. 

\section{Decomposition of torsion and curvature}

On a manifold $M$ with a particular structure group $G$ we can break the
spaces of tensors on $M$ of a particular kind into those taking values
in irreducible subspaces under $G$. For example, on a Riemannian
manifold (of dimension at least 4) the Singer--Thorpe Theorem
\cite{bib:SingThor} says that the curvature of the Levi-Civita
connection breaks into 3 pieces under the orthogonal group, the scalar
curvature, the traceless Ricci tensor and the Weyl curvature. On
oriented 4-manifolds there is a further decomposition of the Weyl
curvature under the special orthogonal group into self-dual and
anti-self-dual parts. We are interested here in the case where $G$ is
the general linear group. Irreducibility is determined by the semisimple
part, and only when counting multiplicities do we need to look at how the
centre acts on summands where highest weight of the semisimple part is
the same.

In the following we suppose that $V$ is a complex vector space and we
handle real vector spaces by replacing them by their complexifications.
Denote by $S^p V$ and $\Lambda^p V$ the $p$-th symmetric and exterior
powers of $V$, respectively, and similarly for the dual space $V^*$. We
denote the symmetric product by juxtaposition and exterior
multiplication by a wedge. The Lie algebra $\gl(V)$ is isomorphic to
$V^*\otimes V$.

As noted in the previous section, torsion tensors are sections of a
bundle whose fibre is associated to the representation of $GL(V)$ on
$\Lambda^2V^*\otimes V$ and curvature tensors the representation  on
$\Lambda^2V^*\otimes \gl(V) = \Lambda^2V^*\otimes V^* \otimes V$.

Let us summarise the highest weight theory that we need to decompose these
representations. If $\dim V = n$ then $\sl(V)$ has rank $n-1$. There are
$n-1$ fundamental representations with highest weights $\omega_i$ and we
number them so that $\omega_i$ is the highest weight of $\Lambda^i V$.
As representations of $SL(V)$, $\Lambda^i V^* = \Lambda^{n-i}V$. The weights
of $V$ are $\{\omega_1, \omega_2 - \omega_1, \omega_3 - \omega_2,\ldots,
\omega_{n-1} - \omega_{n-2}, - \omega_{n-1}\}$ and the weights of $V^*$
are $\{- \omega_1, \omega_1 - \omega_2, \omega_2 - \omega_3,\ldots,
\omega_{n-2} - \omega_{n-1}, \omega_{n-1}\}$.

We denote by $V(m_1,\ldots,m_{n-1})$ the irreducible representation of
$SL(V)$ whose highest weight is $m_1\omega_1+ \dots +
m_{n-1}\omega_{n-1}$. Thus $V= V(1,0,\ldots,0)$, $V^* = V(0,\ldots,0,1)$
$\Lambda^2V^* = V(0,\ldots,0,1,0)$ and so on.

The highest weights occurring in the decomposition of the tensor product
of two irreducibles can be found amongst the highest weight of one
factor plus an arbitrary weight of the other. When the second factor has
weights of multiplicity $1$ then it follows easily from Littelmann's
method \cite{bib:Littelmann} that the set of irreducible components is in
bijection with the dominant elements of this set and these are the
highest weights of the components. Fortunately, both $V$ and $V^*$ have
all weights of multiplicity $1$.

\subsection{Torsion}\label{subsec-tor}
Torsion tensors live in $\Lambda^2V^*\otimes V$ and so the
highest weights of the irreducible factors of this space will be the
dominant elements of the set $\omega_{n-2} + \{\omega_1, \omega_2 -
\omega_1, \omega_3 - \omega_2,\ldots, \omega_{n-1} - \omega_{n-2}, -
\omega_{n-1}\}$. Dominant weights have all coefficients non-negative and
so, by inspection, these are just $\{\omega_{n-2} + \omega_1,
\omega_{n-1}\}$. The second of these is the highest weight of $V^*$.
Thus we have $\Lambda^2V^*\otimes V = \T_1 \oplus
\T_2$  with $\T_1 = V(1,0,\ldots,0,1,0)$ and
$\T_2 = V(0,\dots,0,1) = V^*$.

We have a map $\Lambda^2V^*\otimes V \to V^*$ given by taking a basis
$e_i$ for $V$ and dual basis $\epsilon^i$ for $V^*$ and setting
\[
\widehat{T}(X) = \sum_i \epsilon^i(T(X,e_i))
\]
which defines an element $\widehat{T}$ for each element $T$ of
$\Lambda^2V^*\otimes V$. Conversely, given $\alpha \in V^*$ we can
obtain $\overline{\alpha} \in \Lambda^2V^*\otimes V$ by setting
\[
\overline{\alpha}(X,Y) = \alpha(X)Y - \alpha(Y)X.
\]
Then
\[
\widehat{\overline{\alpha}}(X) = \sum_i \epsilon^i(\alpha(X)e_i -
\alpha(e_i)X) = (n-1)\alpha(X).
\]
It follows that the component of $T$ in $\T_2$ is
\[
T_2(X,Y) = \frac{1}{n-1} \sum_i\epsilon^i(T(X,e_i)) Y - \epsilon^i(T(Y,e_i))X
\]
and the component of $T$ in $\T_1$ is $T_2=T-T_1$. $T= T_1 + T_2$
is then the decomposition of the torsion into irreducible components under
$GL(V)$.

\subsection{Curvature}

We decompose curvature by first looking at $\Lambda^2V^*\otimes V^*$.
This has highest weights the dominant elements amongst $\omega_{n-2} + \{-
\omega_1, \omega_1 - \omega_2, \omega_2 - \omega_3,\ldots, \omega_{n-2}
- \omega_{n-1}, \omega_{n-1}\}$ which are $\{\omega_{n-3}, \omega_{n-2} +
\omega_{n-1}\}$. The first of these corresponds with $\Lambda^3 V^*
=V(0,\ldots,0,1,0,0) $ and the second is a representation we call
$\B_0 = V(0,\ldots,0,1,1)$.

Curvatures live in the space $\Lambda^2V^*\otimes V^* \otimes V$ which
thus has a partial decomposition $\Lambda^3V^* \otimes V \oplus
\B$ where $\B = \B_0 \otimes V$ To decompose
this further we can apply the method again. The first piece decomposes
as the irreducibles which have a highest weight the dominant elements of
$\omega_{n-3} + \{\omega_1, \omega_2 - \omega_1, \omega_3 -
\omega_2,\ldots, \omega_{n-1} - \omega_{n-2}, - \omega_{n-1}\}$, and
these are $\{\omega_{n-3} + \omega_1,\omega_{n-2}\}$. $\B$
decomposes as the irreducibles which have a highest weight the dominant
elements of $\omega_{n-2} + \omega_{n-1} + \{\omega_1, \omega_2 -
\omega_1, \omega_3 - \omega_2,\ldots, \omega_{n-1} - \omega_{n-2}, -
\omega_{n-1}\}$, and these are $\{\omega_{n-2} + \omega_{n-1} +\omega_1,
2\omega_{n-1},\omega_{n-2}\}$. The last two are highest weights of
$S^2 V^*$ and $\Lambda^2 V^*$ respectively.

\subsection{Zero Torsion}
We look at the case where the torsion vanishes. In this case
(\ref{Bianchi1}) implies that the curvature satisfies the first Bianchi
Identity
\[
R^\nabla(X,Y)Z + R^\nabla(Y,Z)X + R^\nabla(Z,X)Y=0.
\]
But, when we have an element $\alpha \in \Lambda^2V^*\otimes V^*$, then
the combination $\frac13 (\alpha(X,Y)Z + \alpha(Y,Z)X + \alpha(Z,X)Y)$
is alternating and is precisely its projection into $\Lambda^3V^*$. Thus
curvatures of torsion zero connections lie in the subspace $\B$
which we call the \textit{Bianchi tensors\/}.

On the space $\B$ the Ricci contraction produces an element
$r^\nabla$ in $V^*\otimes V^* = \Lambda^2 V^* \oplus S^2 V^*$. We want
to go in the reverse direction and build a Bianchi tensor from an
element of $V^*\otimes V^*$.  For this, we use the bracket of
section~\ref{sec:bundle-lie-algebras}: view $Q\in V^*\otimes V^*$ as
a $V^*$-valued $1$-form and the identity map $\Id$ as a $V$-valued
$1$-form then $[Q\wedge\Id]\in\Lambda^2 V^*\otimes\End(V)$ defined by
\[
[Q\wedge\Id](X,Y)=[Q(X),Y]-[Q(Y),X]
\]
takes values in $\B$.  Moreover, the Ricci contraction of
$[Q\wedge\Id]$ is 
\[
-(n+1)Q_--(n-1)Q_+
\]
where $Q_+,Q_-$ are the symmetric and skew parts of $Q$.  Thus the
decomposition of $\B$ into irreducibles reads
\[
\B=\W\oplus[S^2V^*\wedge\Id]\oplus[\Lambda^2V^*\wedge\Id],
\]
with $\W$ the kernel of the Ricci contraction, and the
corresponding decomposition of a curvature tensor $R^\nabla$ is
\[
R^\nabla=W^\nabla-\frac1{n-1}[r^\nabla_+\wedge\Id]-
\frac1{n+1}[r^\nabla_-\wedge\Id].
\]
By analogy with the Riemannian case, we call $W^\nabla$ so defined,
the \emph{Weyl component} of the curvature.

\subsection{Projective Invariance of the Weyl Tensor}

The projective class of a torsion free connection $\nabla$ consists
of all (necessarily) torsion-free connections of the form
\[
\nabla^\alpha=\nabla-[\alpha,\Id],
\]
for a $1$-form $\alpha$.  Thus
\[
\nabla^\alpha_XY=\nabla_XY+\alpha(X)Y+\alpha(Y)X.
\]
A standard computation gives
\begin{align*}
R^{\nabla^\alpha}&=R^{\nabla}-d^\nabla[\alpha,\Id]+
\frac12[[\alpha,\Id]\wedge[\alpha,\Id]]\\
&=R^{\nabla}-[d^\nabla\alpha\wedge\Id]+
\frac12[[\alpha,\Id]\wedge[\alpha,\Id]]
\end{align*}
since $d^\nabla\Id=T^\nabla=0$.  As for the zero-order term, since
$TM$ is abelian, $[\Id\wedge\Id]=0$ and the Jacobi identity (for the
superalgebra of Lie algebra valued forms) then gives
\begin{align*}
[\Id\wedge[\alpha\wedge\Id]]&=[[\Id\wedge\alpha]\wedge\Id]\\
&=-[\Id\wedge[\alpha\wedge\Id]]
\end{align*}
so that $[\Id\wedge[\alpha\wedge\Id]]=0$.  Bracketing this last with
$\alpha$ gives
\[
0=[\alpha,[\Id\wedge[\alpha\wedge\Id]]]=
[[\alpha,\Id]\wedge[\alpha,\Id]]+[\Id\wedge[\alpha,[\alpha\wedge\Id]]]
\]
and we conclude that
\[
R^{\nabla^\alpha}=R^{\nabla}-
[(d^\nabla\alpha+\frac12[\alpha,[\alpha\wedge\Id]])\wedge\Id].
\]
Thus $R^{\nabla^\alpha}-R^{\nabla}$ lies entirely in $[T^*M\otimes
T^*M\wedge\Id]$ and in particular
\[
W^{\nabla^\alpha}=W^\nabla.
\]
Thus the Weyl curvature is a projective invariant.

\section{Application to Twistor Theory}

Additional details on the structure of twistor spaces needed to prove
the results in this section can be found in \cite{bib:OBR,bib:Raw}.

If $M$ is a manifold we denote by $J(M)$ the bundle over $M$ whose fibre
at $x \in M$ consists of all endomorphisms $j$ of the tangent space at
$x$ with $j^2 = -1$, and we let $\pi \colon J(M) \to M$ be the bundle
projection. For this to make sense, the dimension $n$ of $M$ must be
even, say $n=2m$. The differential $d\pi$ of $\pi$ is a surjective
bundle morphism from $TJ(M)$ to the pull-back $E$ of $TM$ to $J(M)$. The
kernel of $d\pi$ is the vertical tangent bundle $\V$. At $j \in
J(M)$ the vertical space $\V_j$ can be identified with
endomorphisms of $E_j$ which anticommute with $j$. $\End(E)$ has a
canonical section $\J$ given by $\J_j=j$.

If $\nabla$ is a connection in $TM$ it induces a pull-back connection
$\pi^*\nabla$  in $E$ and the covariant differential
$\pi^*\nabla\J$ is an $\End(E)$-valued $1$-form on $J(M)$ whose
values anticommute with $\J$. In fact, if we identify
$\V$ with such endomorphisms, then $\pi^*\nabla\J 
\colon TJ(M) \to \V$ is surjective, so the kernel is a subbundle
$\H$ of $TJ(M)$ which is mapped isomorphically onto $E$ by $d\pi$.
$\J$ gives $E$ a complex structure and left multiplication by 
$\J$ gives $\V$ a complex structure. There is then a unique
almost complex structure $J^\nabla$ on $J(M)$ such that
\[
d\pi (J^\nabla X) = \J d\pi (X) \qquad\mathrm{and}\qquad
\pi^*\nabla_{J^\nabla X}\J = \J \pi^*\nabla_{X}\J.
\]

The $(1,0)$ tangent spaces of $J^\nabla$ on $J(M)$ are spanned by
vectors of the form $(J^\nabla + i)X$ and the $(0,1)$ tangent spaces by
$(J^\nabla - i)X$. Thus the $(1,0)$ forms are spanned by components of
$d\pi\circ (J^\nabla + i) = (\J  + i)d\pi$ and 
$X \mapsto \pi^*\nabla_{(J^\nabla + i)X}\J  = 
(\J  + i)\pi^*\nabla_{X}\J$.

If we change the connection from $\nabla$ to $\nabla' = \nabla + A$ then
$\pi^*\nabla' = \pi^*\nabla + \pi^*A$. Thus
\[
\pi^*\nabla'_{}\J  = \pi^*\nabla_{}\J  + 
[\pi^*A, \J ]
\]
and
\[
(\J  + i)\pi^*\nabla'_{}\J  = 
(\J  + i)\pi^*\nabla_{}\J  + (\J  + i)[\pi^*A, 
\J ].
\]
$\nabla'$ will define the same almost complex structure as $\nabla$
provided $(\J  + i)[\pi^*A, \J]$ is a $(1,0)$-form.
This will be the case if and only if
\[
(\J  + i)[\pi^*A((J^\nabla - i)X), \J ] = 
(\J  + i)A((\J  - i)d\pi X)(\J  - i)
= 0.
\]
As an endomorphism-valued $1$-form, $(j + i)A_x((j - i) X)(j - i)$ is
the projection of $A_x$ into the $3i$ eigenspace of $j$ on 
$T^*_xM \otimes T^*_xM \otimes T_xM$

\begin{prop}\label{prop-same-j}
Two connections $\nabla$ and $\nabla' = \nabla + A$ define the same
almost complex structure on $J(M)$ if and only if no irreducible
component of $A$ has values in an irreducible subspace of $V^* \otimes
V^* \otimes V$ where $j_0$ has a $3i$ eigenvalue.
\end{prop}
\begin{pf}
Each $j$ is obtained from one fixed $j_0$ by conjugation by an element
$g$ of $GL(V)$. So a statement about all $j$ is equivalent to a
statement about a single $j_0$ provided we apply it to whole irreducible
components with respect to $GL(V)$ at a time. The result now follows
from the preceding calculation.
\end{pf}

\begin{thm}\label{thm-proj-inv}
The almost complex structure $J^\nabla$ on $J(M)$ defined by a
connection $\nabla$ on $M$ only depends on the projective class of
$\nabla$. If $\nabla$ and $\nabla'$ both have zero torsion then
$J^\nabla = J^{\nabla'}$ if and only if $\nabla$ and $\nabla'$ are
projectively equivalent.
\end{thm}

\begin{pf}
If $\nabla$ and $\nabla'$ are projectively equivalent then $A \in [T^*M
\wedge \Id]$ and so has values in the irreducible components where only
$\pm i$ eigenvalues occur, so $J^\nabla = J^{\nabla'}$ by Proposition
\ref{prop-same-j}.

What remains is to show that for torsion zero connections $J^\nabla =
J^{\nabla'}$ implies that $\nabla$ and $\nabla'$ are projectively
equivalent. But when $\nabla$ and $\nabla' = \nabla + A$ both have zero
torsion, then $A$ has values in the bundle associated to $S^2V^* \otimes
V$. The highest weights of irreducible components will be the dominant
elements of the set $2\omega_{n-1} + \{\omega_1, \omega_2 - \omega_1,
\omega_3 - \omega_2,\ldots, \omega_{n-1} - \omega_{n-2}, -
\omega_{n-1}\}$ and these are $\{\omega_{n-1}, 2\omega_{n-1} +
\omega_1\}$. The first is $V^*$ which has no $3i$ eigenvalue and the
second is $V(1,0,\ldots,0,2)$ which therefore does have a $3i$
eigenvalue. It follows from Proposition \ref{prop-same-j} that $J^\nabla
= J^{\nabla'}$ implies that $A$ has values in the irreducible component
of $S^2V^* \otimes V$ corresponding to $V^*$. This is embedded via
$\alpha \in V^* \mapsto \alpha(X)Y+\alpha(Y)X$. Hence $\nabla$ and
$\nabla'$ are projectively equivalent.
\end{pf}

When is the almost complex structure $J^\nabla$ integrable?
This question was considered in various cases in \cite{bib:OBR}. The
condition discovered there is the analogue for the Nijenhuis tensor of
$J^\nabla$ of Proposition \ref{prop-same-j} and tells us that the
torsion and curvature of $\nabla$ should have irreducible components
only in subrepresentations of $\Lambda^2 V^*\otimes V$ and $\Lambda^2
V^*\otimes \gl(V)$ where $j_0$, viewed as an element of the Lie algebra
$\gl(V)$ acting in these representations has no $3i$ or
$4i$ eigenvalue, respectively. Let us look at the torsion
condition first, as it leads to a simplification of the curvature case.

\begin{thm}\label{thm-torzero}
If $J^\nabla$ is integrable then there is another connection
$\nabla'$ defining the same almost complex structure on $J(M)$ and with
zero torsion.
\end{thm}

\begin{pf}
The eigenvalues of $j_0$ are imaginary and so we really need to work
with the complexification of $V$, but to keep the notation simple we
shall still refer to this as $V$. The eigenvalues of $j_0$ on $V$ and
$V^*$ are $\pm i$, and on a $k$-fold tensor product of these we
get eigenvalues which are a sum of $k$ of these. Thus possible
eigenvalues on torsion tensors are $\pm3i,\pm i$, and the
$3i$ does occur on the whole space as it is the value of the
highest weight on $j_0$. It cannot occur on $\T_2$ since this
is isomorphic to $V^*$ and so the $3i$ eigenvalue occurs on
$\T_1$. By the result of \cite{bib:OBR} it follows that when
$J^\nabla$ is integrable, $T^\nabla$ must lie in the
$\T_2$ subspace. If we set
\[
\alpha(X) = \frac{1}{n-1} \sum_i \epsilon^i(T(X,e_i)),
\]
in the notation of Section \ref{subsec-tor}, then we have
\[
T^\nabla(X,Y) = \alpha(X)Y - \alpha(Y)X.
\]
Consider
\[
\nabla'_XY = \nabla_XY - \frac12(\alpha(X)Y - \alpha(Y)X).
\]
$\nabla'$ clearly has torsion zero and differs from $\nabla$ by
an endomorphism-valued $1$-form with no component having eigenvalue
$3i$. So by Proposition \ref{prop-same-j} $J^\nabla=J^{\nabla'}$.
\end{pf}

Theorem \ref{thm-torzero} allows us to assume, without loss of generality,
that the connection $\nabla$ defining the almost complex structure has
torsion zero, which we do from now on.

\begin{thm}\label{thm-intgrbl}
Let $\nabla$ be a torsion-free connection then $J^\nabla$ is
integrable if and only if $W^\nabla = 0$.
\end{thm}

\begin{pf}
The  curvature argument is similar to the torsion case with the
appropriate changes of representation and eigenvalue. Integrability
forces any irreducible component of the curvature to vanish if the
$4i$ eigenvalue of $j_0$ occurs on the corresponding irreducible
component of $\B$. It is easy to see that (i) the $4i$
eigenvalue does occur on $\B$, that (ii) the only eigenvalues
on $V^*\otimes V^*$ are $0$ and $\pm2i$ and hence that the
$4i$ eigenvalue must actually occur on the space of Weyl
tensors. Thus the integrability condition for $J^\nabla$ is
that the Weyl component, $W^\nabla$, must vanish.
\end{pf}

\begin{rem}
When $\dim M=2$, any $W^\nabla=0$ since then the space of curvature
tensors is just $[T^*M\otimes T^* M\wedge\Id]$.  Thus, in this case,
any $J^\nabla$ is integrable.
\end{rem}

\section{Projectively Flat Connections}

\begin{defn}
Say that a connection $\nabla$ on a manifold
$M$ is \textit{projectively flat} or \textit{of Ricci type} if 
$W^\nabla=0$.
\end{defn}

We are going to show that when $\dim M=n\geq 3$, a projectively flat
connection $\nabla$ induces a local diffeomorphism between $M$ and
$\R P^n$ which intertwines the projective class of $\nabla$ with the
projective class of the Levi-Civita connection on $\R P^n$.

We begin by describing the relevant geometry of $\R P^n$.

Denote by $V$ the trivial bundle $\R P^n\times\R^{n+1}$ and let
$\Lambda\to\R P^n$ be the tautological subbundle of $V$ whose fibre at
$\ell\in\R P^n$ is $\ell\subset V$.  The flat connection $d$ on $V$
induces a canonical isomorphism
$\beta:T\rpn\to\Hom(\Lambda,V/\Lambda)$ such that, for $\sigma$ a
section of $\Lambda$,
\[
\beta(X)\sigma=d_X\sigma\mod\Lambda.
\]
There is a dual isomorphism (also called $\beta$) from $T^*\rpn$ to
$\Hom(V/\Lambda,\Lambda)$ determined by
\[
\beta(\alpha)\beta(X)\sigma=\alpha(X)\sigma
\]
or, equivalently,
\[
\beta(\alpha)d_X\sigma=\alpha(X)\sigma.
\]

The connections in the projective class we wish to consider are in
bijective correspondence with complements to $\Lambda$ in $V$.
Indeed, if $U$ is such a complement so that $V=\Lambda\oplus U$, then
$d$ followed by projection along $\Lambda$ or $U$ gives connections $\nabla$
on $\Lambda$ and $U\cong V/\Lambda$ and so connections on
$\Hom(\Lambda,V/\Lambda)$ and so, via $\beta$, on $T\rpn$.

To compute the curvature and torsion of such a $\nabla$, first note
that $\beta$ induces an isomorphism $T\rpn\Lambda\cong U$ via
\[
X\otimes\sigma\mapsto\beta(X)\sigma\in V/\Lambda\cong U
\]
so that we have a connection-preserving isomorphism
$V\cong\Lambda\oplus T\rpn\Lambda$ with respect to which the flat
derivative decomposes as
\[
d=
\begin{pmatrix}
\nabla&Q\\\Id&\nabla
\end{pmatrix}
\]
for some $Q$ a $1$-form with values in $T^*\rpn$.  Explicitly,
\[
d_X(\sigma,Y\otimes\tau)=(\nabla_X\sigma+Q(X,Y)\tau,
\nabla_X(Y\otimes\tau)+X\otimes\sigma).
\]
We compute the curvature of $d$:
\[
0=R^d=
\begin{pmatrix}
R^\nabla+[Q\wedge\Id]&d^\nabla Q\\
d^\nabla\Id&R^\nabla+[Q\wedge\Id]
\end{pmatrix}.
\]
In particular, $d^\nabla\Id=T^\nabla=0$ so that $\nabla$ is
torsion-free.  Further, $R^\nabla+[Q\wedge\Id]=0$ on both $\Lambda$
and $T\rpn\Lambda$ and so on $T\rpn$ also (we are in the situation of
section~\ref{sec:bundle-lie-algebras} so there is no ambiguity in the
definition of our brackets).  In particular, $R^\nabla$ has no
component in $\W$ so that $\nabla$ is projectively flat.

To see how $\nabla$ varies with $U$, we need an explicit formula for
$\nabla$.  For this, let $\pi:V\to\Lambda$ be the projection along
$U$.  Then
\begin{align*}
\beta(\nabla_XY)\sigma&=(\nabla_X\beta(Y))\sigma\\
&=\nabla_X(\beta(Y)\sigma)-\beta(Y)\nabla_X\sigma\\
&=d_X((1-\pi)d_Y\sigma)-d_Y(\pi d_X\sigma)\mod\Lambda\\
&=d_Xd_Y\sigma-d_X(\pi d_Y\sigma)-d_Y(\pi d_X\sigma)\mod\Lambda.
\end{align*}
The projection along any other complement differs from $\pi$ by a
section of $\Hom(V/\Lambda,\Lambda)$ and so is of the form
$\pi-\beta(\alpha)$ for some $1$-form $\alpha$.  If $\nabla^\alpha$
is the corresponding connection then we have
\begin{align*}
\beta(\nabla^\alpha_X Y-\nabla_XY)\sigma&=
d_X(\beta(\alpha)d_Y\sigma)+d_Y(\beta(\alpha)d_X\sigma)\mod\Lambda\\
&=\beta(X)(\alpha(Y)\sigma)+\beta(X)(\alpha(Y)\sigma)\\
&=\beta([X,\alpha]Y)\sigma.
\end{align*}
Thus $\nabla^\alpha=\nabla-[\alpha,\Id]$.

We have therefore constructed a projective class of projectively
flat torsion-free connections on $\rpn$.

\begin{rem}
Equip $V$ with a flat metric and take $U=\Lambda^\perp$.  The metric
induced on $\Hom(\Lambda,U)$ gives, via $\beta$, an
$\mathrm{SO}(n+1)$-invariant Riemannian structure on $\rpn$ which is
that for which $\rpn$ is a Riemannian symmetric space.  Moreover, the
corresponding connection $\nabla$, being torsion-free and clearly
metric, is the Levi-Civita connection for the symmetric metric.
\end{rem}

Suppose now that $M$ is an $n$-manifold, $n\geq2$, with torsion-free
connection $\nabla$ and let $\Lambda$ be a bundle of
$-n/(n+1)$-densities.  Then $\nabla$ induces a connection on
$\Lambda$, also called $\nabla$, with curvature $F^\nabla$ given by
\[
F^\nabla=-\frac1{n+1}\Tr R^\nabla=-\frac1{n+1}s^\nabla=
-\frac2{n+1}r^\nabla_-.
\]
Moreover, set
\[
Q^\nabla=\frac{1}{n-1}r^\nabla_++\frac{1}{n+1}r^\nabla_-
\]
so that, on $TM$, $R^\nabla+[Q^\nabla\wedge\Id]=W^\nabla$.

Set $V=\Lambda\oplus TM\Lambda$ and, as before, define a connection
$\ca$ on $V$ by
\[
\ca=
\begin{pmatrix}
\nabla&Q^\nabla\\\Id&\nabla
\end{pmatrix}.
\]
Then
\[
R^\ca=
\begin{pmatrix}
F^\nabla+[Q^\nabla\wedge\Id]&d^\nabla Q^\nabla\\
d^\nabla\Id&
(R^\nabla+[Q^\nabla\wedge\Id])\otimes\Id+\Id \otimes 
(F^\nabla+[Q^\nabla\wedge\Id])
\end{pmatrix}.
\]
Once more, $d^\nabla\Id=T^\nabla=0$ while, for $\sigma$ a section of
$\Lambda$,
\begin{align*}
[Q^\nabla\wedge\Id](X,Y)\sigma&= Q^\nabla_X(Y\otimes\sigma)-
Q^\nabla_Y(X\otimes\sigma)\\
&=(Q^\nabla(X,Y)-Q^\nabla(Y,X))\sigma
\end{align*}
so that $[Q^\nabla\wedge\Id]=2Q^\nabla_-=-F^\nabla$ on $\Lambda$. 
Therefore,
\[
R^\ca=
\begin{pmatrix}
0&d^\nabla Q^\nabla\\0&W^\nabla\otimes\Id
\end{pmatrix}.
\]
\begin{rem}
Although it is not at first apparent, our construction of $(V,\ca)$
depends only on the projective class of $\nabla$.  From a more
invariant view-point, $V^*$ is the bundle $J_1(\Lambda^*)$ of
$1$-jets of sections of $\Lambda$ and $\ca$ is the normal Cartan
connection thereon, c.f. \cite{bib:Baston,bib:Cartan}.
\end{rem}

\begin{prop}
For $n=\dim M\geq3$, $\ca$ is flat if and only if $W^\nabla=0$ 
\end{prop}
\begin{pf}
The only issue is to show that when $W^\nabla=0$ then $d^\nabla Q^\nabla$
vanishes also.  Now the second Bianchi identity together with
$T^\nabla=0$ gives
\[
d^\nabla W^\nabla=d^\nabla [Q^\nabla\wedge\Id]=[d^\nabla Q^\nabla\wedge\Id].
\]
Now, for any $2$-form $\omega$ with values in $T^*M$, one readily
computes that
\[
\alpha^i([\omega\wedge\Id](X,Y,X_i)Z)=
-((n-2)\omega_+ + (n+1)\omega_-)(X,Y,Z)
\]
where $\omega_+$ and $\omega_-$ are, respectively, the components of
$\omega$ in $\B_0$ and $\Lambda^3T^*M$.  Thus, when $n\geq3$
and $W^\nabla$ vanishes, both $(d^\nabla Q^\nabla)_\pm$ vanish whence
$d^\nabla Q^\nabla=0$ and $\ca$ is flat.
\end{pf}

\begin{rem}
When $n=2$, any connection is projectively flat so that
$C^\nabla=d^\nabla Q^\nabla$ is the only obstruction to flatness of
$\ca$.  In this case, $C^\nabla$ is an invariant of the projective
class of $\nabla$: it is the \emph{projective covariant} of
Veblen--Thomas \cite{bib:Veblen-Thomas}.  This is the projective
analogue of the Cotton--York tensor of conformal geometry.
\end{rem}

Suppose now that $\ca$ is flat.  On a simply connected open subset
$\Omega$ of $M$, we can trivialise the pair $(V,\ca)$: that is, there
is a bundle isomorphism $\Phi:V\vert_\Omega\cong M\times\R^{n+1}$ such
that $\Phi\circ\ca=d\circ\Phi$.  We now have a map
$\phi:\Omega\to\rpn$ given by $\phi(x)=\Phi\Lambda_x\subset\R^{n+1}$
for which
\[
\phi^{-1}\Lambda_{\rpn}=\Phi\Lambda.
\]
Set $U=\Phi(TM\Lambda)$: a complement to $\Phi\Lambda$.  Then we may
view $\phi^{*}\beta$ as taking values in
$\Hom(\Phi\Lambda,U)\cong\phi^{-1}T\rpn$.  We have
\[
d_X(\Phi\sigma)=\Phi(\ca_X\sigma)
\]
and taking the component in $U$ yields
\[
\phi^*\beta(X)\Phi\sigma=\Phi(X\otimes\sigma).
\]
We conclude that $\phi^*\beta$ is an isomorphism so that $\phi$ is a
local diffeomorphism.  Moreover, $\phi^*\beta$ and hence $d\phi$
intertwines $\nabla$ with the connection on $\rpn$ induced by $U$.

To summarise:
\begin{thm}[\cite{bib:Weyl}]\label{thm-projflat}
Let $M$ be an $n$-dimensional manifold with torsion-free connection
$\nabla$.  Suppose that either $n\geq3$ and $W^\nabla=0$ or $n=2$ and
$C^\nabla=0$.  Then there are local affine diffeomorphisms
between $M$ and $\rpn$ equipped with a connection in the projective
class of the symmetric connection.
\end{thm}

In particular, in this case, the projective class of $\nabla$
contains a locally symmetric connection and the twistor space of $M$
is locally biholomorphic to that of $\rpn$.

\section*{Acknowledgements}

The first author thanks David Calderbank for helpful conversations on
Cartan Geometry.

The second author thanks Roger Carter for explaining the use of Littelmann's
method \cite{bib:Carter,bib:Littelmann} for decomposing the tensor
product of irreducible representations, and Christian Duval and the
CPT, Marseille-Luminy for its hospitality during part of this work.

We thank Izu Vaisman for bringing his work on transversal twistor spaces
of foliations \cite{bib:Vaisman3} to our attention which has a
substantial overlap with the results of Section 4. The work for our
present paper was completed in November 2000, but publication delayed
due to the pressures of other work.

We also thank Thomas Mettler for bringing the paper of Veblen--Thomas
\cite{bib:Veblen-Thomas} to our attention.


\begin{thebibliography}{10}

\bibitem{bib:AHS}
M.~F. Atiyah, N.~J. Hitchin, and I.~M. Singer, \emph{Self-duality in
  four-dimensional {R}iemannian geometry}, Proc. Roy. Soc. London Ser. A
  \textbf{362} (1978), no.~1711, 425--461. \MR{MR506229 (80d:53023)}

\bibitem{bib:Baston}
R.~J. Baston, \emph{Almost {H}ermitian symmetric manifolds. {I}. {L}ocal
  twistor theory}, Duke Math. J. \textbf{63} (1991), no.~1, 81--112.
  \MR{MR1106939 (93d:53064)}

\bibitem{bib:BBO}
L.~B{\'e}rard-Bergery and T.~Ochiai, \emph{On some generalizations of the
  construction of twistor spaces}, Global Riemannian geometry (Durham, 1983),
  Ellis Horwood Ser. Math. Appl., Horwood, Chichester, 1984, pp.~52--59.
  \MR{MR757205 (86h:53028)}

\bibitem{bib:Cartan}
E.~Cartan, \emph{Sur les vari\'et\'es \`a connexion projective}, Bull. Soc.
  Math. France \textbf{52} (1924), 205--241. \MR{MR1504846}

\bibitem{bib:Carter}
R.~W. Carter, \emph{Representations of simple {L}ie algebras: modern variations
  on a classical theme}, Algebraic groups and their representations (Cambridge,
  1997), NATO Adv. Sci. Inst. Ser. C Math. Phys. Sci., vol. 517, Kluwer Acad.
  Publ., Dordrecht, 1998, pp.~151--173. \MR{MR1670769 (2000b:17014)}

\bibitem{bib:Dubois-Violette}
Michel Dubois-Violette, \emph{Structures complexes au-dessus des vari\'et\'es,
  applications}, Mathematics and physics (Paris, 1979/1982), Progr. Math.,
  vol.~37, Birkh\"auser Boston, Boston, MA, 1983, pp.~1--42. \MR{MR728412
  (85h:53053a)}

\bibitem{bib:Eisenhart}
Luther~Pfahler Eisenhart, \emph{Riemannian geometry}, Princeton Landmarks in
  Mathematics, Princeton University Press, Princeton, NJ, 1997, Eighth
  printing, Princeton Paperbacks. \MR{MR1487892 (98h:53001)}

\bibitem{bib:Littelmann}
Peter Littelmann, \emph{Paths and root operators in representation theory},
  Ann. of Math. (2) \textbf{142} (1995), no.~3, 499--525. \MR{MR1356780
  (96m:17011)}

\bibitem{bib:OBR}
N.~R. O'Brian and J.~H. Rawnsley, \emph{Twistor spaces}, Ann. Global Anal.
  Geom. \textbf{3} (1985), no.~1, 29--58. \MR{MR812312 (87d:32054)}

\bibitem{bib:Raw}
John~H. Rawnsley, \emph{{$f$}-structures, {$f$}-twistor spaces and harmonic
  maps}, Geometry seminar ``Luigi Bianchi'' II---1984, Lecture Notes in Math.,
  vol. 1164, Springer, Berlin, 1985, pp.~85--159. \MR{MR829229 (87h:58048)}

\bibitem{bib:Salamon}
S.~M. Salamon, \emph{Quaternionic structures and twistor spaces}, Global
  Riemannian geometry (Durham, 1983), Ellis Horwood Ser. Math. Appl., Horwood,
  Chichester, 1984, pp.~65--74. \MR{MR757207}

\bibitem{bib:SingThor}
I.~M. Singer and J.~A. Thorpe, \emph{The curvature of {$4$}-dimensional
  {E}instein spaces}, Global Analysis (Papers in Honor of K. Kodaira), Univ.
  Tokyo Press, Tokyo, 1969, pp.~355--365. \MR{MR0256303 (41 \#959)}

\bibitem{bib:TricVanh}
Franco Tricerri and Lieven Vanhecke, \emph{Curvature tensors on almost
  {H}ermitian manifolds}, Trans. Amer. Math. Soc. \textbf{267} (1981), no.~2,
  365--397. \MR{MR626479 (82j:53071)}

\bibitem{bib:Vaisman1}
Izu Vaisman, \emph{Symplectic curvature tensors}, Monatsh. Math. \textbf{100}
  (1985), no.~4, 299--327. \MR{MR814206 (87d:53077)}

\bibitem{bib:Vaisman2}
\bysame, \emph{Variations on the theme of twistor spaces}, Balkan J. Geom.
  Appl. \textbf{3} (1998), no.~2, 135--156. \MR{MR1746886 (2001a:53076)}

\bibitem{bib:Vaisman3}
\bysame, \emph{Transversal twistor spaces of foliations}, Ann. Global Anal.
  Geom. \textbf{19} (2001), no.~3, 209--234. \MR{MR1828080 (2002d:53061)}

\bibitem{bib:Veblen-Thomas}
Oswald Veblen and Joseph~Miller Thomas, \emph{Projective invariants of affine
  geometry of paths}, Ann. of Math. (2) \textbf{27} (1926), no.~3, 279--296.
  \MR{MR1502733}

\bibitem{bib:Weyl}
Hermann Weyl, \emph{Zur infinitesimalgeometrie: Einordnung der projektiven und
  der konformen auffassung}, G\"ott. Nachr. (1921), 99--112.

\end{thebibliography}
\providecommand{\bysame}{\leavevmode\hbox to3em{\hrulefill}\thinspace}
\providecommand{\MR}{\relax\ifhmode\unskip\space\fi MR }
\providecommand{\MRhref}[2]{%
  \href{http://www.ams.org/mathscinet-getitem?mr=#1}{#2}
}
\providecommand{\href}[2]{#2}

\end{document}